%% `supleg.tex': Hamiltonian versus Lagrangian
%%             formulations of supermechanics.
%%  (27 july 1996)   
 
\magnification\magstep1
\overfullrule=0pt
\scrollmode

\font\eightrm=cmr8                %% for the abstract
\font\tensmc=cmcsc10              %% for the authors

\font\tenbbb=msbm10 \font\sevenbbb=msbm7 \font\fivebbb=msbm5
\newfam\bbbfam
\textfont\bbbfam=\tenbbb \scriptfont\bbbfam=\sevenbbb
  \scriptscriptfont\bbbfam=\fivebbb
\def\Bbb{\fam\bbbfam}      %% Blackboard bold (math mode only)
 
\font\tengoth=eufm10 \font\sevengoth=eufm7 \font\fivegoth=eufm5
\newfam\gothfam
\textfont\gothfam=\tengoth \scriptfont\gothfam=\sevengoth
  \scriptscriptfont\gothfam=\fivegoth
\def\goth{\fam\gothfam}    %% Euler Fraktur (math mode only)
 
\def\opname#1{\mathop{\rm#1}\nolimits} %% operator name
 
\def\a{\alpha}                    %% abbreviation for  \alpha
\def\A{{\cal A}}                  %% Sheaf
\def \AEU{\A_E{\cal (U)}}         %%  abbreviation of {\A_E{\cal (U)}}
\def \AU{{\cal A(U)}}             %%  abbreviation of {\cal A(U)} 
\def\abs#1{\vert#1\vert}          %% substitute for  ,z,
         %% anticommutator {X,Y}
              %% adjoint action
\def\B{{\cal B}}                  %% abbreviation of {\cal B}
                  %%%* (the original \b )
\def\b{\beta}                     %%* abbreviation for  \beta
                                  %%% (Prefer \b to \be -- Joe)
         %% closure operator
     %% <w,A,z>
 %% <w,z>
  %% doubled bracket <<w,z>>
\def \BU{{\cal B(U)}}             %%  abbreviation of {\cal B(U)}
\def\bw{\bigwedge\nolimits}       %% big exterior product
\def\C{{\cal C}}                  %% canonical section in T(TM)
\def\cite#1{\lbrack{\bf#1}\rbrack} %% reference citation
    %%* equivalence class  ,X|
 %% commutator  ,X,Y|
                  %% set of Sode's
\def\de{\delta}                   %% abbreviation for  \delta
\def\De{\Delta}                   %% abbreviation for  \Delta
\def\Der{\opname{Der}}            %%* derivations
\def\dim{\opname{dim}}            %% dimension
\def\E{{\cal E}}                  %% abbreviation for  {\cal E}
            %% endomorphism ring
\def\ENE{{\cal N}}                %% abbreviation for  {\cal N}
\def\eq#1{{\rm(#1)}}              %% print equation number
\def\FL{{\cal FL}}                %% Legendre transformation
                 %% abbreviation for  \varphi
\def\Ga{\Gamma}                   %% abbreviation for  \Gamma
\def\ga{\gamma}                   %% abbreviation for  \gamma
               %% general linear algebra
\def\GL{\opname{GL}}              %% general linear group
                 %% Heisenberg superalgebra
\def\harr#1#2{\smash{\mathop{\hbox to .5in{\rightarrowfill}}
     \limits^{\scriptstyle#1}_{\scriptstyle#2}}}
                    %%horizontal arrangement in square a diagram
\def\harrt#1#2{\smash{\mathop{\overline{\hbox to .5in{\hrulefill}}}
     \limits^{\scriptstyle#1}_{\scriptstyle#2}}}
                %%horizontal arrangement in triangular a diagram
          %% abbreviation for  \widehat
              %% fraction  1/2
\def\Hom{\opname{Hom}}            %% endomorphism ring
\def\ID{\relax{1\kern-.24 em\rm l}}  %% identity in a superalgebra
\def\I{{\cal I}}                  %% superideal
              %%* identity morphism
         %% inclution of T^3\M
\def\Im{\opname{Im}}              %% image of an operator
 %% implication arrow
\def\ker{\opname{ker}}            %% kernel of an operator
                  %% Lie derivative
                  %% exterior algebra
                  %% abbreviation for  \lambda
\def\M{{\cal M}}                  %% Supermanifold
               %% omit next symbol
\def\Maps{\opname{Maps}}          %% morphisms of graded manifolds
                      %% abbreviation for \eta
                  %% nonnegative integers
         %% substitute for  \|F\|
\def\o{\omega}                    %% abbreviation for \omega
\def\om{\omega}                    %% abbreviation for \omega
\def\O{\Omega}                    %% abbreviation for \Omega
                %% orthogonal Lie algebra
       %% abbreviation for \widetilde
             %%* orthosymplectic Lie salg
\def\ox{\otimes}                  %%* tensor product  (Joe)
\def\o+{\oplus}                   %% direct sum
\def\pd#1#2{{{\partial#1}\over{\partial#2}}}  %% partial derivative
               %% general linear superalgebra
             %%* \pi as binary operator
\def\Proof{\noindent{\sl Proof\/}}   %% proof label
                  %% sheaf of a superdomain
\def\R{{\Bbb R}}                  %% real numbers
                  %% skewsymmetric operators
\def\sepword#1{\qquad\hbox{#1}\quad} %%* well-spaced words
\def\set#1{\{\,#1\,\}}            %% set notation
  %% tiny fraction  1/2
\def\si{\sigma}                   %% abbreviation for  \sigma
\def\Si{\Sigma}                   %% abbreviation for  \Sigma
        %% sign function
 %% strong operator limit
\def\smc{\tensmc}                 %% small capitals
          %% vector space span
              %% symplectic Lie algebra
\def\sopd#1#2#3{{{{\partial^2{#1}}\over{{\partial{#2}\partial{#3}}}}}}    
                                  %% second order partial derivative
\def\srule#1#2{(-1)^{\vert#1\vert\,\vert#2\vert}} 
                                  %%* sign rule in superalgebras
\def\sruleiii#1(#2#3){(-1)^{\abs{#1}\,(\abs{#2}+\abs{#3})}}
                                  %%* sign rule with 3 terms
             %% supertrace
 %%* str of operator
          %% support of
               %% trace
\def\th{\theta}                   %% abbreviation of \theta
\def\Th{\Theta}                   %% abbreviation of \Theta
\def\T{{\bf T}}                   %% canonical vertor field along
                         %% the proyection of the tangent bundle
\def\Tau{{\cal T}}                %% abbreviation for {\cal T}
                %% kth-prolongation of \T
            %% (k-1)th-prolongation of \T
            %% abbreviation for \T^{(1)}
            %% abbreviation for \T^{(2)}
 %% small fraction  1/2
   %% small fraction  i/2
 %% small fraction  1/4
 %%* trace of operator
\def\U{{\cal U}}                  %% open set of M
                      %% symmetric product
                  %% open set of N
\def\vphi{\varphi}                %% abbreviation for \varphi
\def\vth{\vartheta}               %% abbreviation of \vartheta
\def\varr#1#2{\llap{$\scriptstyle #1$}\left\downarrow
     \vcenter to .5in{}\right.\rlap{$\scriptstyle #2$}}
            %%vertical arrangement for a square diagram with a down arrow.

            %%vertical arrangement for a square diagram with an up arrow.          
                    %% exterior product
         %% abbreviation; usage: \wrt~
\def\x{\times}                    %% cartesian product
\def\X{{\goth X}}                 %% vector fields
                %% second prolongation of X
        %%  coordinates of 2nd prolongation of X
                %% kth-prolongation of X
            %% (k-1)th-prolongation of X
                %% first prolongation of X
        %%  coordinates of 1st prolongation of X
\def\z{\zeta}                     %% abbreviation for \zeta
                  %% integers
\def\1{\'{\i}}                    %% Spanish accented i
\def\3{\sharp}                    %% abbreviation for # symbol
\def\7{\dagger}                   %% abbreviation for + symbol
\def\.{\cdot}                     %% anonymous variable
\def\:{\colon}                    %% colon in  f: A -> B
\def\<#1,#2>{\langle#1\mathbin\vert#2\rangle} %% (= `\braket')

\def\qed{\allowbreak\qquad\null
         \nobreak\hfill\square}   %% end-proof marker at right
\def\square{\vbox{\hrule
             \hbox{\vrule height 5.2pt
              \hskip 5.2pt
              \vrule}\hrule}}     %% wizard's blank box, for show

\def\today{\number\day\space      %%* date, English style
           \ifcase\month\or
            January\or February\or
            March\or April\or May\or
            June\or July\or August\or
            September\or October\or
            November\or December\fi
            \number\year}

\newcount\secnum
\outer\def\beginsection#1. #2\par{%  %% Redefine the Plain Tex
      \vskip 0pt plus.3\vsize        %%  `\beginsection'
      \penalty -250                  %%  to set the section
      \vskip 0pt plus-.3\vsize       %%  number for automatic
      \bigskip \vskip\parskip        %%  equation numbering.
      \global\secnum=#1%
      \global\equationnum=0
      \message{#1. #2}%
      \leftline{\bf#1. #2}\nobreak
      \smallskip\noindent}
 
\outer\def\subsection#1.#2. #3\par{%
      \ifnum#2=1 \smallskip       %% Subsection divisions
       \else \bigskip \fi         %%*  (spacing corrected)
      \message{Sec #1.#2.}%
      \leftline{\it#1.#2. #3}\nobreak
      \smallskip\noindent}
 
\def\declare#1. #2\par{\medskip   %%* for definitions,
          \noindent{\bf#1.}\rm    %%  remarks, examples, etc.
          \enspace\ignorespaces
          #2\par\smallskip}
 
         %%* items with separation

\newcount\equationnum             %% We combine ``Fisher's rule''
\def\ecnum{\the\secnum.%          %%  with cross-referencing;
           \the\equationnum}      %%  say `\eqlabel\cs'
\everydisplay={\global\advance    %%  AFTER any display
               \equationnum by 1  %%  and refer to it thereafter
               \leftnumdisplay}   %%  by `(\cs)' or `\eq\cs'.
\def\leftnumdisplay#1{\let\nxxt#1\maybealign}
\def\maybealign{\ifx\nxxt\leqalignno
                \else\leftnumeqn\fi\nxxt}
\def\leftnumeqn#1$${#1\leqno(\ecnum)$$}  %% automatic numbering
 
\def\eqlabel#1{\xdef#1{\ecnum}}   %% assign label to latest eqn
 
\def\refno#1. #2\par{\smallskip   %% format references
          \item{\lbrack#1\rbrack}
                #2\par}
 
\long\def\suspend#1\resume{}      %% hide some input paragraphs

\def\Abraham{1}
\def\Batchelor{2}
\def\Sancheziii{3}
\def\Hector{4}
\def\Hectoriii{5}
\def\Pepino{6}
\def\Pepinii{7}
\def\Pepiniv{8}
\def\Gawedzki{9}
\def\Hilli{10}
\def\Hillii{11}
\def\Ibort{12}
\def\Kostant{13}
\def\Leites{14}
\def\Manin{15}
\def\MonterdeSanchez{16}
\def\Monterde{17}
\def\Sanchez{18}
\def\Sanchezi{19}
\def\Sanchezii{20}
\def\Tennison{21}
\def\Yano{22}

%%% Document begins here:  
\centerline{\bf Hamiltonian versus Lagrangian formulations of
supermechanics.}
 
\bigskip 
\centerline{\smc Jos\'e F. Cari\~nena} 
\smallskip 
\centerline{\it Departamento de F\'{\i}sica Te\'orica, Universidad de 
Zaragoza, 50009 Zaragoza, Spain.}  
\medskip \centerline{\smc H\'ector Figueroa} 
\smallskip 
\centerline{\it Escuela\ de Matem\'atica, Universidad de Costa                
Rica, San Jos\'e, Costa Rica.}
 
\bigskip {\narrower\eightrm\baselineskip=9.5pt           
%% (Abstract)
 
\noindent We take advantage of different generalizations of the tangent
manifold to the context of graded manifolds, together with the notion of 
super section along a morphism of graded manifolds, to obtain intrinsic 
definitions of the main objects in supermechanics such as, the vertical 
endomorphism, the canonical and the Cartan's graded forms, 
the total time derivative operator 
and the super--Legendre transformation. In this way, we obtain a 
correspondence 
between the Lagrangian and the Hamiltonian formulations of
supermechanics.\par} 

\bigskip\noindent{\it Keywords\/}: Tangent superbundle, Supervector field 
along a morphism, Super--Legendre transformation. 
 
\noindent{\it 1991 MSC numbers\/}: Primary: 58A50, 58C50. Secondary: 70H33 

\noindent{\it PACS numbers\/}:  03.20+i, 02.90+p, 11.30.pb.
\vfill
\vfill
\vfill
To appear in J. Phys. A: Math. Gen.
\vfil\eject 
\noindent {\bf 1. Introduction.} 
\medskip
 The  idea of considering classical systems that incorporate
commuting and anticommuting variables to study dynamical systems
dealing with bosonic and fermionic degrees of freedom, in particular
supermechanics, has been in the air for some time now.  Moreover, it has 
proved to be quite useful, not only in physics but also in
mathematics.  Nevertheless, a careful study of the  geometric
foundations of supermechanics was not taken very seriously, or at
least people did not pay
 the necessary attention, until quite
 recently~\cite{\Ibort}, in spite of the general tendency to geometrize 
 physics.  One of the reasons for this is that although
 the general consent is that the proper setting is the theory of
 supermanifolds, there is no general agreement, for instance, as to
 what the velocity phase space of the system should be, since
 there are several different possibilities to generalize the
 concept of tangent bundle in the context of graded manifolds.   
 One of the central points in~\cite{\Ibort} was the introduction
of the tangent supermanifold, which proved to be the right arena
to develop Lagrangian supermechanics, since it allowed an
intrinsic theory. However, some of the central objects,
although well defined, were not defined in an intrinsic way. 
 Perhaps the main 
drawback of the tangent supermanifold is that it is not a bundle.  
To overcome this, we enlarge this tangent supermanifold
by considering the tangent superbundle as introduced by S\'anchez--Valenzuela 
in~\cite{\Sanchez}, which unfortunately is a
little too big, as its dimension is $(2m+n,2n+m)$ if the dimension of
the starting graded manifold (the superconfiguration space) is $(m,n)$, 
but that has the big advantage of allowing a geometric
interpretation of a supervector field as a section of a superbundle
in much the same way as in non--graded geometry.  We shall show in
this paper the convenience of getting a compromise between both
concepts: we shall introduce the objects using the tangent
superbundle approach, but thereafter
 we shall read the results in the tangent supermanifold (identified 
 as a subsupermanifold of the tangent superbundle).  It will be 
 shown how the tangent superbundle
 structure is the appropriate framework for an intrinsic definition
 of objects such as the total time derivative operator, the
 vertical superendomorphism, the Cartan 1--form and, fundamentally,
 the Legendre transformation, which will allow us to establish  a
 connection between the Lagrangian and the Hamiltonian formalisms of 
 supermechanics, similar 
to the one in classical mechanics.  

 In the geometrical approach to 
classical mechanics, the infinitesimal transformations arising in the 
traditional approach are described by the flow of vector fields, 
which can be considered either as sections of the tangent bundle, or as 
derivations of the commutative algebra $C^\infty(M)$ of differentiable 
functions. The generalization of the concept of a flow of
a supervector field is not an easy task \cite{\MonterdeSanchez},
but the corresponding idea of vector field translates easily to the 
framework of supermechanics.
It was shown in~\cite{\Pepino} that in order to incorporate 
non--point transformations in velocity phase space, it is necessary 
to introduce the concept of supervector field along a map.  Moreover, the
use  of such concept
and its generalizations, sections of a vector bundle along a map,
has proved to be very useful for a better understanding of many
aspects of classical mechanics~\cite{\Pepinii,\Pepiniv}.
 What we want  to remark is that  in the
transition to the supermechanics setting these concepts are even
more necessary because of the inconvenience of working with points
in graded geometry.  Therefore, in the process of constructing a 
geometrical approach to supermechanics, including fermionic degrees
of freedom, one of the first concepts to be introduced is that of
section along a morphism of supermanifolds. 

   The organization 
of the paper is as follows: in Section 2 we 
describe the tangent superbundle, in particular we give a ``Batchelor's 
description" of it, and discuss its relation to
the tangent supermanifold as defined by Ibort and Mar\1n in \cite{\Ibort}; 
it is shown that supervector fields can then be seen as
geometric sections of the tangent superbundle. In Section 3, we introduce 
the notion of a section along a morphism of graded 
manifolds, and represent
supervector fields along a morphism as geometric sections along the 
morphism of the tangent superbundle and, as a particular
example, we give an intrinsic definition of the total time derivative
operator that was used in~\cite{\Hector} to obtain a version of Noether's
theorem in supermechanics, and that plays an important role in the
geometry of the tangent superbundle, and thereby in the Lagrangian
Formalism of supermechanics.  

Section 4 is devoted to the study of
graded forms along a morphism of graded manifolds.  Furthermore,
we study the canonical graded 1--form $\Th_0$ on the supercotangent  
manifold, as well as the degeneracy of the graded form $\O_0= -d\Th_0$. 
Finally, Section 5  
is concerned with the vertical superendomorphism,
which is necessary to introduce the Cartan 1--form corresponding to
a Lagrangian superfunction, and also with the super Lengendre transformation. 
  Finally using the machinery here developed,
we establish a relationship between the Lagrangian and the 
Hamiltonian formulations of supermechanics.

\beginsection 2. The tangent superbundle and the tangent supermanifold.

\subsection 2.1. Basic notation.   

At the heart of the graded manifold 
theory is the idea of equiping a supervector space $V=  V_0\oplus V_1$ 
with the structure of a graded manifold; the natural way of doing 
this~\cite{\Kostant,\Leites} is to consider the so called 
affine supermanifold:
$$
S(V) := \bigl
(V_0,C^\infty(V_0) \ox \bw(V_1^*)\bigr).
$$  
Nevertheless, this has some 
drawbacks from the categorical
point of view~\cite{\Sanchezi}, and, in the context 
of supervector bundles, S\'anchez--Valenzuela 
realized that it is more appropriate to use, 
instead of the affine supermanifold, what he 
coined~\cite{\Sanchez,\Sancheziii} the 
supermanifoldification of $V$:
$$
V_S := S(V\oplus \Pi V),
$$
where $\Pi$ is the change of parity 
functor~\cite{\Leites,\Manin}, hence 
$(\Pi V) = (\Pi V)_0 \oplus (\Pi V)_1$, where
$$
(\Pi V)_i = V_{i+1}  \qquad i = 0,1.
$$  
The sheaf 
$C^\infty(V_0) \ox \bw(V_1^*)$ will be 
denoted by $\A_{m,n}$ whenever $\dim V_0 = m$, 
$\dim V_1 = n$, and $\R^{m|n}$ will denote the 
graded manifold $\R^{m|n} = (\R^m,\A_{m,n})$.  
On the other hand, we shall always  consider, 
on $\R^{m+n|m+n} = (\R^m \oplus \R^n)_S$, the
following supercoordinates: if 
$\set{e_i,r_\a: i=1,\dots,m, \quad \a=1,\dots,n}$ 
is a graded basis for 
$\R^m \oplus \R^n$ (so $\abs{e_i} = 0$ 
and $\abs{r_\a} = 1$) and 
$\set{t^i,\vth^\a}$ is the corresponding 
dual basis, then the 
set $\set{t^i,\pi\vth^\a;\vth^\a,\pi t^i}$ 
gives a supercoordinate system in $\R^{m+n|m+n}$.  
Here $\pi$ is the natural morphism between $V$ and 
$\Pi V$.

\subsection 2.2. The supertangent bundle.

  Our first goal is to describe the relation 
between the
supertangent manifold as defined 
in~\cite{\Ibort,\Hector} and
the supertangent bundle introduced 
by S\'anchez--Valenzuela
in~\cite{\Sanchez}.

  If $\M = (M,\A_M)$ is a graded manifold 
of dimension $(m,n)$,
its supertangent bundle is defined 
via the one--to--one
correspondence between equivalence classes of 
locally free
sheaves of $\A_M$--modules over $\M$ of rank $(r,s)$,
and equivalence classes of 
supervector bundles over $\M$ of rank $(r,s)$, 
considered 
as a natural generalization of the standard 
definition 
of vector bundles; namely, as the quadruplets
$\set{(E,\A_E),\Pi,(M,\A_M),V_S}$ such that
$\Pi \: (E,\A_E) \to (M,\A_M)$ is a 
submersion of graded 
manifolds, $V$ is a real $(r,s)$--dimensional 
supervector
space and every $q\in M$ lies in a 
coordinate neighbourhood
$\U \subseteq M$ for which an isomorphism 
$\Psi_\U$ exists 
making the following diagram commutative:
$$
\matrix{
\Bigr(\pi^{-1}(\U),\A_E\bigr(\pi^{-1}(\U)\bigl)\Bigl)
&\harr{\displaystyle \Psi_\U}{}  
&\bigr(\U,\A_M(\U)\bigl) \x V_S    \cr
\varr{\displaystyle\Pi}{}  
&& \varr{}{\displaystyle P_1}  \cr
\bigr(\U,\A_M(\U)\bigl)  
&\harrt{}{}  &\bigr(\U,\A_M(\U)\bigl).\cr} 
$$

In fact, the supertangent bundle is defined 
precisely as the
supervector bundle of rank $(m,n) = \dim \M$ 
that corresponds
to the sheaf of $\A_M$--modules $\Der \A$.

   As the superbundle $(E,\A_E)$ is locally 
isomorphic to a
graded manifold of the form 
$\bigl(\U,\AU\bigr) \x V_S$, we shall take
advantage of this fact to describe the local 
supercoordinates
of $(E,\A_E)$.  Thus, if $\set{q^i,\th^\a}$, 
$i=1,\dots, m$,
$\a= 1,\dots, n$, are  local supercoordinates on 
$\U\subseteq M$,
and $\set{t^j,\pi\vth^\b,\vth^\b,\pi v^j}$, 
$j=1,\dots, r$,
 $\b= 1,\dots, s$, are the local supercoordinates of 
$V_S = \R^{m+n|m+n}$ described previously, then
$\set{p^*_1q^i,p^*_2t^j,p^*_2\pi\vth^\b,
p^*_1\th^\a,p^*_2\vth^\b,p^*_2\pi t^j}$, where 
$P_k =(p_k,p^*_k)$ 
is the natural projection of $\bigl(\U,\AU\bigr) \x V_S$ 
onto the $k$--th factor, is a set of local 
supercoordinates on $\bigl(\U,\AU\bigr) \x V_S$,
hence the image of this set under the 
morphism of superalgebras
$\psi^*$ will be a set of local supercoordinates 
for $(E,\A_E)$
on $\pi^{-1}(\U)$, which, abusing of the 
notation, we shall denote by 
$\set{q^i,v^j,\pi\z^\b,\th^\a,\z^\b,\pi v^j}$.

\declare Remark 2.1.
 We also want to point out that the superideal 
$\I$, locally generated by the superfunctions 
$\set{\pi v^j,\pi \z^\b}$
($1\leq j\leq r$ and $1\leq \b\leq s$), defines 
a subsupermanifold
of $(E,\A_E)$ of dimension $(m+r,n+s)$.  
Similarly, the superideal
$\I'$ locally generated by the superfunctions 
$\set{v^j,\z^\b}$
defines another subsupermanifold of $(E,\A_E)$ 
of dimension $(m+s,n+r)$.

\subsection 2.3. Simple graded manifolds.

  Next we want to describe the supertangent bundle 
$ST\M := (STM, ST\A)$ in a more concise way.  With this
in mind, we shall first make some comments on the
Batchelor--Gawedzki structural theorem.
Let $\pi\: E\to M$ be a vector bundle of rank $n$, and
$\bigwedge E$ its exterior algebra vector bundle (i.e.
the vector bundle over $M$ whose fiber on a point
$q\in M$ is the vector space $\bigwedge E_q$). The sheaf
of sections $\Ga(\bigwedge E)$ can be considered, in the
obvious way, as a sheaf of supercommutative superalgebras
over $M$.  Moreover, $\bigl(M,\Ga(\bigwedge E)\bigr)$
is a graded manifold.  Indeed, if
$\set{(\U_k,\phi_k)}$ is an
atlas of $M$ such that $\pi^{-1}(\U_k)$
trivialize $E$, then we have diffeomorphisms 
$\phi_k\: \U_k \to U_k\subseteq\R^m$ and
$\psi_k\: \pi^{-1}(\U_k) \to \U_k\x\R^n$
such that $pr_1\circ\psi_k = \pi|_{\U_k}$.  
Consider the superdomain
$\bigl(U_k,\A_{m,n}(U_k)\bigr)$
and let $\set{u^i_k,\xi^\a_k}$ ($i=1,\ldots,m$
and $\a=1,\dots,n$), be supercoordinates on it.
 Now, if 
$\th^\a_k\:\U_k\to\pi^{-1}(\U_k)$ 
is the local section of $\bw E$ defined by 
$\th^\a_k(u)= \psi^{-1}_k(u,e_\a)$, where 
$\set{e_1,\dots,e_n}$ denotes the canonical basis 
of $\R^n$,
it is clear that the morphism $\Phi_k\: 
\Bigl(\U_k,\Ga\big(\bw\pi^{-1}(\U_k)\bigr)\Bigr)
\to \bigl(U_k,\A_{m,n}(U_k)\bigr)$
defined by the assignments
$$
u^i_k \mapsto q^i_k:= \pi_i\circ\phi_k,
\sepword{and}
\xi^\a_k \mapsto \th^\a_k,
$$
where $\pi_i\: \R^m \to \R$ is the 
projection onto the $i$--th factor,
is a chart, in the sense of graded manifolds, for 
$\bigl(M,\Ga(\bw E)\bigr)$.  Moreover, it is easy to
check that, if  $\U_{kl} := \U_k \bigcap\U_l
\ne\emptyset$  
the transition function of this graded
manifold
$$
\Phi_{kl}\:\Bigl(\phi_l(\U_{kl}),
\A_{m,n}\bigl(\phi_l(\U_{kl})\bigr)\Bigr) 
\longrightarrow \Bigl(\phi_k(\U_{kl}),
\A_{m,n}\bigl(\phi_k(\U_{kl})\bigr)\Bigr),
$$
is given by the relations
$$
\leqalignno{
\phi_{kl}^*(u^i_k) &= \phi^i_{kl} 
& (\ecnum\rm a) \cr
\phi_{kl}^*(\xi^\a_k) &=
(\psi_{lk})_{\b\a}\,\xi^\b_l, 
& (\ecnum\rm b) \cr}
\eqlabel\cero
$$
where  $\phi_{kl} = \phi_k \circ \phi_l^{-1}$
denotes the change of coordinates in $M$,  
$\psi_{kl} = \psi_k \circ \psi_l^{-1}$ 
is the transition function of the vector bundle 
$\pi\: E\to M$ over $\U_{kl}$, and
$(\psi_{kl})_{\a\b}$ is the  matrix associated
to  $\psi_{kl}$. We refer to this kind of  graded
manifolds as simple graded manifolds.

 Simple graded manifolds are more than just a nice
example of graded manifolds. Indeed, it is not hard
to obtain a fiber bundle out of a graded manifold.
  Let  $\set{\U_j}$ be an open cover of
$M$ such that on each $\U_j$ one has local charts
of $\M$, say 
$\Phi_j\: \bigl(\U_j,\A_M(\U_j)\bigr) \to
\bigl(U_j,\A_{m,n}(U_j)\bigr)$ and let 
$\set{u^i_j,\xi^\a_j}$  ($i=1,...,m$ and
$\a=1,...,n$) be supercoordinates on
$\bigl(U_j,\A_{m,n}(U_j)\bigr)$. If the
transition morphisms are given by the relations
$$
\leqalignno{
\phi_{jk}^*(u^i_j)  &= (\phi_{jk})^i_0(u) + 
(\phi_{jk})^i_{\a\b}(u)\,\xi^\a_k\xi^\b_k + \cdots
& (\ecnum\rm a)   \eqlabel\ocho   \cr
\phi_{jk}^*(\xi^\a_j) &=
(\varphi_{jk})^\a_\b(u)\,\xi^\b_k + 
(\varphi_{jk})^\a_{\b\ga\de}(u)\,
\xi^\b_k\xi^\ga_k\xi^\de_k  + \cdots,
& (\ecnum\rm b)    \cr}
\eqlabel\muno
$$
then, from the cocycle relations of the
$\Phi_{jk}$'s it follows that the
matrices $(\varphi_{jk})_{\a\b}$
satisfy, on each point of 
$\phi_k(\U_j \bigcap \U_k \bigcap\U_l)$, the
cocycle relations
$$
\varphi_{jk} \circ \varphi_{kl} = \varphi_{jl}.
$$
Thus the functions 
$\tilde\vphi_{jk} \:\U_j \bigcap \U_k \to
\GL(n,\R)$, defined by $\tilde\vphi_{jk}(q) =
\Bigl(\vphi_{jk}\bigl(\phi_k(q)\bigr)
\Bigr)_{\a\b}$, give rise to a vector bundle
$E \to M$.
  Now, if we also asume that the $\U_j$'s are
such that the $\pi^{-1}(\U_j)$'s
trivialize $E \to M$, then by our previous
argument we have a local chart $\Psi_j\: 
\Bigl(\U_j,\Ga\bigl(\bw\pi^{-1}(\U_j)\bigr)
\Bigr) \to \bigl(U_j,\A_{m,n}(U_j)\bigr)$ of
$\bigl(M,\Ga(\bw E)\bigr)$.  Moreover,
$\psi_j^* \circ (\phi_j^*)^{-1}$ is an
isomorphism from the superalgebra $\A(\U_j)$ 
into the superalgebra
$\Ga\bigl(\bw\pi^{-1}(\U_j)\bigr)$. Thus,
the graded manifolds
$(M,\A_M)$ and $\bigl(M,\Ga(\bw E)\bigr)$ are
locally isomorphic.  Surprisingly enough, these 
graded
manifolds are globally isomorphic, although not
in a canonical way, a fact known as the
structural theorem of
Batchelor~\cite{\Batchelor} and 
Gawedzki~\cite{\Gawedzki}.

\declare Remark 2.2.
What we want to emphasize is that, from
\eq\cero, an explicit way to construct the
so called structural bundle $E\to M$ is to
use the functions $\vphi_{jk}$, the first
term of the second equation of \eq\muno,
as the transition functions of the dual
bundle $E^*$.

\subsection 2.4. The underlying manifold of
the supertangent bundle.

In order to describe the tangent superbundle
$ST\M := (STM, ST\A)$ we shall follow the
general construction of a supervector bundle out
of a sheaf of $\A_M$--modules 
given in~\cite{\Sanchez} applied to the sheaf
of supervector fields $\Der \A$.
Let $\U$ be an open subset of $M$ such that 
$\bigr(\U,\A_M(\U)\bigl)$ is isomorphic
to a superdomain;
if $X = \sum_{i=1}^m X^i\partial_{q^i}
+ \sum_{\a=1}^n \chi^\a
\partial_{\th^\a}$ is a supervector field
 in $\Der \A(\U)$,  then the map
$$
g_\U \: X \longmapsto (X^1,\cdots,X^m,
\chi^1,\cdots,\chi^n)
$$
defines an isomorphism between the sheaves
of $\A_M$--modules $\A(\U)^m\o+ \A(\U)^n$
and $\Der\AU$.   
Moreover, if $\bigr(\U_1,\A_M(\U_1)\bigl)$ and
$\bigr(\U_2,\A_M(\U_2)\bigl)$
are two of such superdomains then the map
 $$
g_{12} = g_1(\U_1\cap\U_2) \circ
 g_2^{-1}(\U_1\cap\U_2) \:\A(\U_1\cap
 \U_2)^m \o+ \A(\U_1\cap\U_2)^n
\to\A(\U_1\cap\U_2)^m \o+ \A(\U_1\cap\U_2)^n,
$$
which basically expresses the change of
supercoordinates of the supervector field $X$,
is an isomorphism of
$\A(\U_1\cap\U_2)$--modules and is explicitly
given by the graded matrix
$$
g_{12} =\pmatrix{A_{12}  &  \Th_{12}   \cr          
 \Ga_{12}  &  D_{12}  \cr}
 =\pmatrix{\pd{q^i_1}{q^j_2}  &  
 \pd{q^i_1}{\th^\b_2}   \cr          
 \pd{\th^\a_1}{q^j_2}
&  \pd{\th^\a_1}{\th^\b_2} \cr} .
\eqlabel\tres
$$  
 
 Since $g_{12}$ is invertible then the matrices
$\tilde A_{12}$ and  $\tilde D_{12}$, obtained
from $A_{12}$ and $D_{12}$, respectively, 
by projecting their entries onto
$C^\infty(\U_1\cap\U_2)$, are also 
invertible~\cite{\Leites}; moreover, since the
$g$'s satisfy the cocycle condition, we also have
$$
\tilde A_{12} \circ \tilde A_{23} = \tilde A_{13}
\sepword{and}
\quad\tilde D_{12}\circ\tilde D_{23} =\tilde D_{13}.
$$   
The conclusion is that the matrices
$$
\tilde g_{12} =\pmatrix{\tilde A_{12}  &  0   \cr 
0  &  \tilde D_{12}  \cr},
$$
give rise to a smooth vector bundle
$\tau \: STM \to M$, which is the 
Whitney sum of the vector bundle determined by
the transition functions 
$\tilde A_{12} = \pd{q^i_1}{q^j_2}$,
which is nothing but the tangent 
bundle of the manifold $M$, and the vector bundle
$\widetilde E \to M$ determined 
by the $\tilde D$'s, which by the Remark 2.2,
is isomorphic to the dual bundle of the
structural bundle of $\M$. Therefore, we have
proved the following proposition:

\proclaim Proposition 2.1.
If $E\to M$ is a vector bundle such that
$(M,\A) \cong \bigl(M,\Ga\bw(E)\bigr)$, then
the underlying manifold of tangent superbundle
of $\M$ is 
$$
STM = TM \o+ E^*.
$$      \qed

\subsection 2.5. The sheaf $ST\A$.

 To complete the description of the tangent
superbundle we should describe 
the sheaf $ST\A$.  This description is done
in terms of the 
matrices \eq\tres\ taking in consideration
the fact that locally $ST\M$ is isomorphic to  
$\bigr(\U,\AU\bigl) \x  \R^{m+n|m+n}$.
  Thus, if $\tau \: TM \o+ E^* \to M$ is the
canonical projection, then, according to
\cite{\Sanchez}, $ST\A$ is constructed using
the superdomains $\Bigl(\tau^{-1}(\U_j), 
ST\A\bigl(\tau^{-1}(\U_j)\bigr)\Bigr)$ and the
superalgebra morphisms defined by the relations
$$
\leqalignno{
q^i_1 = \hat\Phi_{12}(q^i_1) 
&= \phi^i_0(q) + \phi^i_{\a\b}(q)\,\th^\a_2\th^\b_2 
+ \cdots
& (\ecnum\rm a) \cr
\th_1^\a = \hat\Phi_{12}(\th_1^\a) &=
\psi^\a_\b(q)\,\th^\b_2 
+\psi^\a_{\b\ga\de}(q)\,\th^\b_2\th^\ga_2\th^\de_2 
+\cdots.
& (\ecnum\rm b) \cr
v_1^i = \hat\Phi_{12}(v_1^i)
&= \sum^m_{j=0} \pd{q_1^i}{q^j_2} v^j_2
   - \sum^n_{\b=0} \pd{q_1^i}{\th^\b_2} \z^\b_2
& (\ecnum\rm c) \cr
&= \biggl(\pd{\phi^i_0}{q^j} + 
     \pd{\phi^i_{\a\b}}{q^j}\,\th^\a_2\th^\b_2 
    + \cdots \biggr){v^j_2} 
   + (2\phi^i_{\a\b}(q)\,\th^\a_2 + \cdots)\z^\b_2 \cr
\pi\z_1^\a = \hat\Phi_{12}(\pi\z_1^\a) 
&= -\sum^m_{j=0} \pd{\th_1^\a}{q^j_2} \pi v^j_2 
    + \sum^n_{\b=0} \pd{\th_1^\a}{\th^\b_2}\pi \z^\b_2
& (\ecnum\rm d) \cr
&= -\biggl(\pd{\psi^\a_\b}{q^j_2}\,\th^\b_2 
    + \cdots \biggr) \pi v^j_2 +
(\psi^\a_\b + 3\psi^\a_{\b\ga\de} \,\th^\ga_2\th^\de_2 
    + \cdots)\pi\z^\b_2    \cr 
\z_1^\a = \hat\Phi_{12}(\z_1^\a) 
&= \sum^m_{j=0} \pd{\th_1^\a}{q^j_2} v^j_2 
    + \sum^n_{\b=0} \pd{\th_1^\a}{\th^\b_2} \z^\b_2 
& (\ecnum\rm e) \cr
&= \biggl(\pd{\psi^\a_\b}{q^j_2}\,\th^\b_2  
+ \cdots \biggr)v^j_2 +   (\psi^\a_\b(q) 
+ 3\psi^\a_{\b\ga\de}(q) \,\th^\ga_2\th^\de_2 
    + \cdots)\z^\b_2  \cr
\pi v_1^i = \hat\Phi_{12}( \pi v_1^i)
&= \sum^m_{j=0} \pd{q_1^i}{q^j_2} \pi v^j_2
   + \sum^n_{\b=0} \pd{q_1^i}{\th^\b_2} \pi\z^\b_2
& (\ecnum\rm f) \cr
&= \biggl(\pd{\phi^i_0}{q^j} + 
     \pd{\phi^i_{\a\b}}{q^j}\,\th^\a_2\th^\b_2 
+ \cdots \biggr)\pi v^j_2 
+ (-2\phi^i_{\a\b}\,\th^\a_2 + \cdots) \pi\z^\b_2, \cr}
\eqlabel\mdos
$$
where $\set{q^i_j,v^i_j,\pi\z^\a_j,\th^\a_j,
\z^\a_j,\pi v^i_j}$ are the supercoordinates
on $\tau^{-1}(\U_j)$ described in Section 2.2.

   Now, according to the Remark 2.2, the
transition functions of the structural bundle
$E' \to STM$ of $(STM,ST\A)$ are obtained from
\eq\mdos; actually, they are the inverse
transpose of the linear functions
$\Psi_{12} \: \tau^{-1}(\U_1) \cap 
\tau^{-1}(\U_2)  \to \GL(2n + m,\R)$ given by
$$
\Psi_{12}(q,v,\pi\z) =
\pmatrix{
\psi^\a_\b(q)  &  0  & 0   \cr
\pd{\psi^\a_\b}{q^i}v^i  &  \psi^\a_\b(q) & 0   \cr
 -2\phi^i_{\a\b}(q)\pi\z^\b & 0  &  
\pd{\phi^i_0}{q^j}   \cr}.
\eqlabel\mtres   
$$
Here $\set{q,v,\pi\z}$ are local coordinates on
$STM$. Nevertheless, by our arguments in Section 2.3
(i.e. the Batchelor--Gawedzki theorem) we may 
assume that $\phi^i_{\a\b}(q) = 0$. Then, the
following proposition follows 
immediately from \eq\mtres:

\proclaim Proposition 2.2.
If $E\to M$ is a vector bundle such that
$(M,\A) \cong \bigl(M,\Ga\bw(E)\bigr)$, then
the structural bundle of $ST\M$ is isomorphic 
to $(TE \o+ TM)^* \to TM \o+ E^*$. \qed

  We point out that, using different arguments,
the tangent supermanifold has also been studied
in \cite{\Monterde}.
 
 Finally, we notice that the subsupermanifold
that corresponds to $ST\M$, according to the
Remark 2.1,
is nothing but the tangent supermanifold
$(TM, T\A)$ introduced by 
Ibort and Mar\1n--Solano in~\cite{\Ibort}.  

\subsection 2.6. Supervector fields as geometric sections.   

  The main reason for considering the
tangent superbundle
$\{(STM,ST\A),\Tau,(M,\A),V_S)\}$, and
supervector bundles 
in general~\cite{\Sanchez}, is that
their geometrical sections
are in a one--to--one correspondence with the
sections of the
corresponding locally free sheaf of graded
$\A$--modules;
in our case, with the sections of the sheaf
$\Der \A$, in other words, with the supervector
fields over $\M$. Following~\cite{\Sanchez} 
we will make this correspondence explicit
in the particular case we are 
interested in.
The central point of this correspondence
is to notice that
both, the geometric sections and the
``algebraic" sections, when restricted to an
appropriate open set, are isomorphic to
$\Maps\Bigl(\bigl(\U,\A(\U)\bigr),V_S\Bigr)$
the morphisms between the graded manifolds
$\bigl(\U,\A(\U)\bigr)$ and $V_S$.  
First of all, we notice that
$$
\Der\A(\U) \cong \A(\U)^m \o+ \A(\U)^n
\cong \Maps\Bigl(\bigl(\U,\A(\U)
\bigr),V_S\Bigr).
\eqlabel\uno
$$  

If $X\in\Der\A(\U)$ is written in local coordinates as 
$X = \sum_{i=1}^m X^i\partial_{q^i}
+ \sum_{\a=1}^n \chi^\a\partial_{\th^\a}$, 
then \eq\uno\ is implemented by the maps
$$
X \longmapsto (X^1,\cdots,X^m,\chi^1,\cdots,\chi^n)
\longmapsto 
\Phi_X,
\eqlabel\dos
$$
where $\Phi_X = (\phi_X,\phi_X^*)
\in \Maps\Bigl(\bigl(\U,\A(\U)\bigr),
V_S\Bigr)$ is the morphism described by,
see~\cite{\Leites}, the morphism of 
superalgebras $\phi_X^* \: \A_{m+n,m+n} \to
\A(\U)$ corresponding to the assignments:
$$
\matrix{t^i &\mapsto X^i_0 \qquad\qquad
&\pi\vth^\a &\mapsto \chi^\a_0,   \cr
\vth^\a &\mapsto \chi^\a_1  \qquad\qquad
&\pi t^i &\mapsto X^i_1,   \cr}
$$
where $X^i_0$ denote the even part of
$X^i \in \A(U)$, and so on.    
On the other hand, if $ST\A(\U)$ is a
short notation for
$ST\A\bigl(\tau^{-1}(\U)\bigr)$ and
$F= (f,f^*)$ is a morphism in 
$\Maps\Bigl(\bigl(\U,\A(\U)\bigr),V_S\Bigr)$,
then $\Si_F \: \bigl(\U,\A(\U)\bigr) \to
\bigl(\tilde\tau^{-1}(\U), ST\A(\U)\bigr)$
will denote the section of the tangent
superbundle described by the morphism
of superalgebras $\si_F^* \: ST\A(\U) \to \A(\U)$
defined by the assignments
$$
\matrix{q^i &\mapsto &q^i  \qquad\qquad
&\th^\a &\mapsto &\th^\a, \cr 
v^i &\mapsto &f^*(t^i)  \qquad\qquad
&\z^\a &\mapsto  &f^*(\vth^\a),   \cr
\pi\z^\a &\mapsto &f^*(\pi\vth^\a)  \qquad\qquad
&\pi v^i &\mapsto &f^*
(\pi t^i).   \cr}
$$
(We remind the reader of our notation concerning
supercoordinates described 
in subsections 2.1 and 2.2).
 It is easy to check that
$$
\Der\A(\U) \cong
\Maps\Bigl(\bigl(\U,\A(\U)\bigr),V_S\Bigr)\cong 
\Ga\Bigl(\bigl(\U,\A(\U)\bigr),\bigl(\tilde\tau^{-1}(\U), 
ST\A(\U)\bigr)\Bigr),
$$
is implemented by the morphisms:
$$
X \longmapsto \Phi_X \longmapsto
\Si_X,
$$
where $\si_X^* \: ST\A(\U) \to 
\A(\U)$ is given by the assignments
$$
\matrix{q^i &\mapsto &q^i  \qquad\qquad
&\th^\a &\mapsto &\th^\a, \cr 
v^i &\mapsto &X^i_0  \qquad\qquad
&\z^\a &\mapsto  &\chi^\a_1,   \cr
\pi\z^\a &\mapsto &\chi^\a_0  \qquad\qquad
&\pi v^i &\mapsto &X^i_1.   \cr}
$$

\beginsection 3. Supervector fields along a morphism.

  Since the information of a graded manifold is
concentrated in the algebraic part, that is in the
sheaf of superalgebras, to carry over the point
constructions of the classical geometry in the
graded context is somewhat difficult; for instance
the notion of a flow of a supervector field is far
from trivial \cite{\MonterdeSanchez,\Hilli,\Hillii}.  To tackle
these problems we introduced in \cite{\Hector}
the notion of a supervector
field along a morphism, which also turned out to be a
useful tool to study (higher order) supermechanics
\cite{\Hectoriii}.  Nevertheless, there they were
defined as some kind of superderivations, and
our goal now is to give to such supervector fields a
geometric description similar to the one in non--graded 
geometry.  It is important to point out that, already in
the non--graded context,  vector fields along a map
simplify several constructions \cite{\Pepino--\Pepiniv}.

\subsection 3.1. Definition.

\declare Definition 3.1.
Let $\Phi =(\phi,\phi^*)\: (N,\B) \to (M,\A)$
be a morphism of graded manifolds;
a homogeneous supervector field along $\Phi$ 
is  a  morphism of sheaves over $M$,
$X\: \A \to \Phi_*\B$ such that for each open 
subset $\U$ of $M$
$$
X(fg) = X(f) \, \phi^*_\U(g)
+ \srule{X}{f} \phi^*_\U(f) \, X(g),
$$
whenever $f\in \A(\U)$ is homogeneous of degree
$\abs{f}$. The sheaf of supervector fields 
along $\Phi$ will be denoted by $\X(\Phi)$.

  If $X$ is a supervector field on $(M,\A)$,
then 
$$
\hat X := \phi^* \circ X \in \X(\Phi),
\eqlabel\zero
$$
is a supervector field along $\Phi$.
  In particular, when $\Phi$ is a regular closed 
imbedding \cite{\Leites}, $\hat X$ is nothing but
the restriction of the supervector field
$X$ to the graded submanifold $\ENE$.
 
  If $Y$ is a supervector field on $(N,\B)$, then
$$
T\phi(Y) := Y \circ \phi^*,
\eqlabel\mone
$$
also belongs to $\X(\Phi)$, and we say that 
$Y$ is proyectable with respect to $\Phi$ if
there exists $X\in\X(\A)$ such that
$$
T\phi(Y)= \hat X.
$$

 $\X(\Phi)$ is a locally free sheaf of 
$\Phi_*\B$--modules over $M$ of rank
$(m,n) = \dim\M$~\cite{\Hector}.
Moreover, if $(q^i,\th^\a)$ ($1\le i\le m$,
$1\le \a\le n$), are local supercoordinates on
$\U\subset\M$, then
$$
\partial_{\hat q^i} := \widehat\partial_{ qi}
\qquad\qquad 
\partial_{\hat\th^\a} := \widehat\partial_{\th^\a},
$$
form a local basis of $\X(\Phi)(\U)$.
In particular, any $X \in\X(\Phi)(\U)$
can be written as
$$
X = \sum_{i=1}^m X^i\partial_{\hat q^i}
+ \sum_{\a=1}^n \chi^\a
\partial_{\hat\th^\a},
\eqlabel\UNO
$$     
where $X^i = X(q^i)$ and
$\chi^\a = X(\th^\a)$ are superfunctions in
$\B\bigl(\phi^{-1}(\U)\bigr)$ 
(denoted, from now on, by $\BU$ for short). 

\subsection 3.2.  Supervector fields along a morphism
as sections along a morphism.

  Geometrical sections of a super vector bundle are
definded as usual:

\declare Definition 3.2.
Let $\Phi \: (N,\B) \to (M,\A)$ be a morphism
of graded manifolds  and let
$\set{(E,\A_E),\Pi,(M,\A_M),V_S}$ be a supervector
bundle over $\M$; a local section of $\E := (E,\A_E)$
along $\Phi$ over an open subset $\U$ of $M$ is a
morphism $\Si= (\si,\si^*)\:\big(\phi^{-1}(\U),\BU\bigr)
\to  \big(\pi^{-1}(\U),\AEU\bigr)$, where again
$\AEU := \A_E\bigl(\pi^{-1}(\U)\bigr)$, 
satisfying the condition
$$
\Phi_\U = \Pi_\U \circ \Si_\U;
\eqlabel\DOS
$$
here the subscript $\U$ means the restriction 
of the morphism to the corresponding
open graded submanifold.  The set 
of such sections will be denoted by
$\Ga_\Phi(\Pi|_\U)$.

It is straightforward to check that the assignment
$$
W \longmapsto \Ga_\Phi(\Pi|_W),
$$
for each open set $W \subseteq \U$, makes
$\Ga_\Phi(\Pi|_\U)$
into a sheaf of $\Phi_*\B$--modules.
Moreover, if $\U$ is a trivialising neighbourhood
of the supervector bundle $\E$, then 
it is not hard to obtain a one--to--one correspondence
between $\Ga_\Phi(\Pi|_\U)$ and $\Maps\Bigl
(\bigl(\phi^{-1}(\U),\B(\U)\bigr),V_S\Bigr)$;
in particular, one concludes that
$\Ga_\Phi(\Pi|_\U)$ is locally free.

\declare Remark 3.1.
In the case when the morphism $\Phi$ is the
projection $\Pi$ of the supervector bundle,
there is a canonical section $\C$, to wit
the identity morphism on $\E$. It turns out
that several relevant objects are defined
using this section.  

 We now turn our attention to the case when the
supervector bundle is the tangent superbundle
$ST\M$, in other words, to supervector fields.
The correspondence between supervector fields
along a morphism and sections along a morphism
of the supervector bundle is carried out along
the same lines as in the case of the usual
supervector fields (see Section 2.6). Thus,
one has
$$
\X(\Phi)(\U) \cong \BU^m \o+ \BU^n
\cong \Maps\Bigl(\bigl(\phi^{-1}(\U),
\BU\bigr),V_S\Bigr).
$$ 
This correspondence is also implemented
by \eq\dos, where now the
superfunctions $X^i$ and $\chi^\a$ are
given by \eq\UNO.  
On the other hand, if $\,\U$ is also a
trivialising neighbourhood 
of the supervector bundle $ST\M$,
as before, one can check that
$$
\Maps\Bigl(\bigl(\phi^{-1}(\U),
\BU\bigr),V_S\Bigr)\cong \Ga_\Phi(\Tau|_\U).
$$ 
The explicit correspondence between a
supervector field $X\in\X(\Phi)(\U)$ 
and a local section along $\Phi$ is
given by
$$
X \longmapsto \Si_X,
$$
where $\si^*_X \: ST\AU \to \BU$ is
defined by the assignments
$$
\matrix{q^i &\mapsto &\phi^*(q^i)
\qquad\qquad  &\th^\a &\mapsto 
&\phi^*(\th^\a),    \cr
v^i &\mapsto &X^i_0  \qquad\qquad  &\z^\a 
&\mapsto  &\chi^\a_1,   \cr\pi\z^\a
&\mapsto &\chi^\a_0  \qquad\qquad
&\pi v^i &\mapsto &X^i_1.   \cr}
$$

\subsection 3.3. The total time derivative operator.

   As in the non--graded context, the geometry of
the tangent supermanifold is concentrated in two
objects: the vertical superendomorphism and the
total time derivative operator.  Moreover this
operator, introduced in \cite{\Hector}, turned 
out to be quite important in the Lagrangian
formalism of supermechanics.  In what follows,
we shall use the previous ideas to provide an
intrinsic definition of the total time derivative
operator.   

\declare Definition 3.3.
The canonical section of the tangent supervector
bundle $(ST\M,\Tau,\M)$ along $\Tau$ described
in the Remark 3.1, will be called the total time
derivative operator and will be denoted by $\T$.

  Since $\T$ is nothing but the identity morphism,
$\T$ corresponds, according to the previous section,
to the superderivation along $\Tau$ given, in terms 
of the standard supercoordinates of $ST\M$, by
$$
\T = \sum_{i=1}^m (v^i+\pi v^i)\partial_{\hat q^i} 
+\sum_{\a=1}^n (\z^\a+\pi\z^\a)\partial_{ \hat\th^\a}.
$$

   As we shall see later on, sometimes it is
convenient to work with the tangent supermanifold $T\M$.
Thus, if $\Phi\:T\M\to ST\M$ is the regular closed
imbedding that defines $T\M$ and that is locally defined
by the relations
$$
\matrix{
q^i &\mapsto & q^i  \qquad\qquad  
&\th^\a &\mapsto & \th^\a,   \cr      
v^i &\mapsto &v^i  \qquad\qquad  
&\z^\a &\mapsto  &\z^\a,   \cr
\pi\z^\a &\mapsto &0  \qquad\qquad 
&\pi v^i &\mapsto &0,   \cr}
$$
then the restriction of $\T$ to $T\M$ would be
the superderivation along the restriction of
$\Tau$ to $T\M$ given by $\phi^*\circ\T$, and 
its local expression would be
$$
\T = \sum_{i=1}^m v^i\partial_{\hat q^i} 
+\sum_{\a=1}^n \z^\a\partial_{\hat\th^\a};
\eqlabel\mfour
$$
where now 
$$
\partial_{\hat q^i}=
\phi^*\circ\tau^*\circ\partial_{q^i}
\sepword{and}
\partial_{ \hat\th^\a}=
\phi^*\circ\tau^*\circ\partial_{\th^\a}.
$$
We shall make no distinction in the notation when
we regard $\T$ as an operator either on $ST\M$ or 
on $T\M$.

\beginsection 4. Graded 1--forms along a morphism 
of supermanifolds.

\subsection 4.1. The cotangent superbundle and 
the cotangent supermanifold.

   The sheaf of graded 1--forms is, by definition, 
the dual sheaf of $\Der\A$, and corresponds, 
according to~\cite{\Sanchez}, to a supervector 
bundle $(ST^*\M,\Pi,\M,V_S)$ that will be called the 
cotangent superbundle of $\M$.  As one might expect,
most of the ideas of the previous sections can be
used with this sheaf of $\A$--modulos, taking in
consideration what happens in the non--graded
context.

   Obviously $\O^1(\A)=\X(\A)^*$ is locally free.
Moreover, if $\,\U$ is an open subset of $M$, and 
$\set{q^i,\th^\a}$ are local supercoordinates on it, 
then $\set{dq^1,\dots,dq^m,-d\th^1,\dots,-d\th^n}$ 
is the basis of the module 
$\O^1\AU = \bigr(\Der\AU\bigl)^*$ dual to the basis
$\set{\partial_{q^i},\partial_{\th^\a}}$ of 
$\X\bigl(\AU\bigr)$.  In particular,
any $\om\in \O^1\AU$ can be written in a unique way, 
in the form 
$$
\om = \sum_{i=1}^m w^i\, dq^i 
+ \sum_{\a=1}^n \om^\a \,d\th^\a,
$$     
where the superfunctions $w^i$ and $\om^\a$ are given by
$$
w^i = \om(\partial_{q^i}) \sepword{and}\quad 
\om^\a = -\om(\partial_{\th^\a}).
\eqlabel\one
$$

  Naturally, one can described the cotangent superbundle
in a similar way as we described the tangent superbundle
in Sections 2.4 and 2.5, but, in analogy with the 
non--graded geometry, using instead the matrices
$(g_{\a\b}^{st})^{-1}$, where $g_{\a\b}$ are the transition 
functions for the tangent superbundle \eq\tres\ and $st$ 
denote the supertranspose matrix.

  The correspondence between the sections of the cotangent 
superbundle $ST^*\M = (ST^*M,ST^*\A)$ and graded 1--forms 
is accomplished using the same ideas as in Section 2.6. Thus, 
if in addition, the open subset $\U$ 
is a trivialising neighbourhood for $ ST^*\M$, such that 
$\bigr(\pi^{-1}(\U), ST^*\AU\bigl)$ is also isomorphic 
to a superdomain, where $\Pi=(\pi,\pi^*)$ is the natural 
projection of $ST^*\M$ on $\M$, and $ST^*\AU$ is a short 
notation for $ST^*\A\bigr(\pi^{-1}(\U)\bigl)$,
then the correspondence
$$
\O^1\AU \cong \Ga\Bigl(\bigl(\U,\A(\U)\bigr),
\bigl(\pi^{-1}(\U), ST^*\A(\U)\bigr)\Bigr),
$$
is implemented by the morphism:
$$
\om \longmapsto  \Si_\om,
$$
where $\si_\om^* \: ST^*\A(\U) \to \A(\U)$ is 
defined by the assignments
$$
\matrix{
q^i &\mapsto &q^i  \qquad\qquad  
&\th^\a &\mapsto &\th^\a,  \cr      
p^i &\mapsto &w^i_0  \qquad\qquad  
&\eta^\a &\mapsto  &\om^\a_1,   \cr
\pi\eta^\a &\mapsto &\om^\a_0  \qquad\qquad 
&\pi p^i &\mapsto &w^i_1.   \cr}
$$
where the $w^i$ and the $\om^\a$ are as in \eq\one\ 
and the subindices 0 or 1 stand for the even or odd 
components.  Once more, we remind the reader of our 
notation concerning local supercoordinates of supervector
bundles.

  In analogy with the tangent superbundle, 
the subsupermanifold
$T^*\M = (T^*M,T^*\A)$ of $ST^*\M$, of dimension $(2m,2n)$,
associated to the superideal $\I^*$ locally generated by
the superfunctions $\set{\pi p^i,\pi\eta^\a}$
will be called the cotangent supermanifold.

\subsection 4.2. Graded forms along a morphism.

\declare Definition 4.1.
Let $\Phi \: (N,\B) \to (M,\A)$ be a morphism of 
graded manifolds; we define $\O^1(\Phi)$, 
the sheaf of graded 1--forms along $\Phi$, as the 
sheaf of $\phi_*\B$--modules dual 
to the sheaf $\X(\Phi)$. In other words,
$$ 
\O^1(\Phi) = \X(\Phi)^* = \Hom(\X(\Phi),\phi_*\B).
$$
In general, $k$-superforms are defined as
$$
\O^k(\Phi) := \bw^k\bigl(\O^1(\Phi)\bigr),
$$
where the wedge product is to be understood in the sense of 
graded algebras.

   Since $\O^1(\Phi)$ is the dual of a locally free 
$\phi_*\B$--modulo, is itself a locally free
$\phi_*\B$--modulo.  Moreover, if $\om$ is a graded 
1--form on $\M$, the restriction of $\om$ to $\ENE$ 
is the graded 1--form along $\Phi$ defined by 
$$ 
\hat\om(\hat X) := \phi^* \circ \om (X)
\qquad \forall X\in\X(\A_M). 
$$
If $(q^i,\th^\a)$ are supercoordinates of $\M$ on $\U$, 
and $d\hat q^i$, $d\hat\th^\a$ are the restrictions of
$dq^i$ and $d\th^\a$ respectively, then
$$
d\hat q^i(\partial_{\hat q^j}) = \de_{ij},\qquad   
d\hat q^i(\partial_{\hat\th^\b}) = 0,\qquad
d\hat\th^\a(\partial_{\hat q^j}) = 0,\qquad   
d\hat\th^\a(\partial_{\hat\th^\b}) = -\de_{\a\b};
$$
hence $\set{d\hat q^i , - d\hat\th^\a}$ is the dual
basis of $\set{\partial_{\hat q^i},\partial_{\hat\th^\a}}$.  
In particular, any graded 1--form $\om$ along 
$\Phi$ can be written locally as
$$
\om = \sum_{i=1}^m w^i \,d\hat q^i 
+ \sum_{\a=1}^n \om^\a \,d\hat\th^\a,
$$ 
where the superfunctions $w^i$ and $\om^\a$ belong to 
$\BU$, and are defined by 
$$
w^i = \om(\partial_{\hat q^i}) \sepword{and}\quad 
\om^\a = -\om(\partial_{\hat\th^\a}).
\eqlabel\mfive
$$

  The equivalent process to \eq\mone\ does not work here;
instead, if $\om$ is a graded 1--form along $\Phi$, then 
$\phi^\sharp\om$ given by
$$
\phi^\sharp\om (Y):= \om \bigl(T\phi(Y)\bigr)
\qquad\forall Y\in\X(\B),
\eqlabel\msix
$$
is a graded 1--form on $\ENE$.  As a matter of 
fact, it is possible to classify the graded
1--forms on $\ENE$ that come from graded 1--forms 
along $\Phi$, when $\Phi$ is a submersion. The
result is that $\O^1(\Phi)$ is isomorphic to the
$\phi_*\B$--modulo of $\Phi$--semibasic 1--forms on
$\ENE$ \cite{\Hector}.

  Naturally, this constructions, together with the
last result, can be generalized to graded $k$--forms.
For instance, if $\om\in\O^k(\Phi)$, then 
$$
\phi^\sharp\om (Y_1,\dots ,Y_k) := 
\om\bigl(T\phi(Y_1),\dots ,T\phi(Y_k)\bigr).
$$

  The important point is that these two processes can be 
combined to give an intrinsic definition of the pull back
of a graded form; something that, to our knowlegde, was
lacking in the graded context. 
 
\declare Definition 4.2.
Let $\Phi \: (N,\B) \to (M,\A)$ be a morphism of 
graded manifolds and let $\mu$ be a graded $k$--form 
on $\M$.  The pull back of $\mu$ by $\Phi$ is the graded 
$k$--form on $\ENE$ given by
$$
\Phi^*(\mu) := \phi^\sharp (\widehat\mu).
$$

   If $\mu$ is the graded 1--form given in local 
supercoordinates by $\mu = \sum_{i=1}^m u^i\, dq^i 
+ \sum_{\a=1}^n \mu^\a \,d\th^\a,$ then 
$$
\widehat\mu = \sum_{i=1}^m \phi^*(u^i)\, d\hat q^i + 
\sum_{\a=1}^n \phi^*(\mu^\a) \,d\hat\th^\a;
$$
on the other hand, if $Y\in\X(\ENE)$ is given in local 
coordinates by $Y = \sum_{j=1}^r Y^j\partial_{ p^j} 
+ \sum_{\b=1}^s \Upsilon^\b\partial_{\eta^\b}$, and 
$\phi^i:= \phi^*(q^i)$ and $\phi^\a := \phi^*(\th^\a)$ 
are the coordinate representation of $\Phi$~\cite{\Leites}, 
then
$$
\leqalignno{
\bigl(\Phi^*\mu\bigr)(Y) 
&= \sum_{ij} \phi^*(u^i)Y^j\, \pd{\phi^i}{p^j} +    
\sum_{i\b} \phi^*(u^i)\Upsilon^\b\, \pd{\phi^i}{\eta^\b} 
& (\ecnum) \cr
& \quad + (-1)^{\abs{Y}} \sum_{j\a} 
\phi^*(\mu^\a) Y^j\, \pd{\phi^\a}{p^j} +
(-1)^{\abs{Y}} \sum_{\a\b} 
\phi^*(\mu^\a) \Upsilon^\b\, \pd{\phi^\a}{\eta^\b};  \cr}
$$
which is the definition given in \cite{\Kostant}.

  The following technical result will be needed later on.

\proclaim Lemma 4.1.
Let $\Phi=(\phi,\phi^*)\: (N,\B) \to (M,\A)$  be a 
diffeomorphism, and $\mu$ a graded $k$--form on $\M$, 
then
$$ 
\phi^{-1*}\bigl(\Phi^*\mu(Y_1,\dots,Y_k)\bigr) =
\mu(\phi^{-1*}\circ Y_1 \circ\phi^*,\dots,
\phi^{-1*}\circ Y_k \circ\phi^*)\bigr)
$$

\Proof:
Since $\Phi$ is a diffeomorphism any supervector 
field on $\ENE$ is projectable with respect to $\Phi$;
hence for each $Y_i$ there exists $X_i\in\X(\A)$ on $\M$ 
such that $Y_i\circ\phi^*= \phi^*\circ X_i$, and one has
$$
\leqalignno{
\bigl(\Phi^*\om\bigr)(Y_1,\dots,Y_k) 
&= \bigl(\phi^\sharp(\widehat\om)\bigr)
(Y_1,\dots,Y_k)
& (\ecnum) \cr
&=\widehat\om(Y_1 \circ\phi^*,\dots,
Y_k \circ\phi^*)  \cr
&= \widehat\om(\phi^*\circ X_1 ,\dots,
\phi^*\circ X_k )   \cr
&= \phi^*\bigl( \om(X_1,\dots,X_k)\bigr),   \cr}
$$
and since $X_i=\phi^{-1*}\circ Y \circ\phi^*$ the
lemma follows.  \qed

\subsection 4.3. The canonical graded forms on 
the cotangent supervector bundle.

  As expected, graded 1--forms along a morphism 
have their geometric counterpart.
  If $\Phi \: \ENE \longrightarrow \M$ is a morphism
of graded manifolds and $\U \subseteq \M$ is an open 
subset such that $\bigl(\U,\AU\bigr)$ is isomorphic
to a superdomain and trivialize the cotanget
superbundle $\Pi\: ST^*\M\to\M$, then the
correspondence
$$
\O^1(\Phi)(\U) \cong \Ga_\Phi(\Pi|_\U),
$$
is carried out using similar arguments as before
and is given by
$$
\om \longmapsto \Si_\om,
$$
where $\si^*_\om \: ST^*\AU \to \BU$ is defined by
the assignments
$$
\matrix{
q^i &\mapsto &\phi^*(q^i)  \qquad\qquad  
&\th^\a &\mapsto &\phi^*(\th^\a),   \cr      
p^i &\mapsto &w^i_0  \qquad\qquad  
&\eta^\a &\mapsto  &\om^\a_1,  \cr
\pi\eta^\a &\mapsto &\om^\a_0  \qquad\qquad 
&\pi p^i &\mapsto &w^i_1.   \cr}
\eqlabel\mseven
$$
 Here $w^i$ and $\om^\a$ are the superfunctions defined 
in \eq\mfive, and the subindices 0 and 1 stand for the even
and odd components, repectively.

  Once again, when $\Phi = \Pi= (\pi,\pi^*)$ is the
canonical proyection of $ST^*\M$ on $\M$, we have,
according to the Remark 3.1, a canonical section
along $\Pi$, which, in view of \eq\mseven,
corresponds to the graded 1--form $\check\Th_0$
in $\O^1(\Phi)$ locally given by 
$$
\check\Th_0=\sum_{i=1}^m (p^i +\pi p^i)\, d\hat q^i + 
\sum_{\a=1}^n (\eta^\a +\pi \eta^\a) \,d\hat\th^\a.
\eqlabel\meight
$$

\declare Definition 4.3.
The graded 1--form $\Pi$--semibasic that corresponds 
to $\check\Th_0\in\O^1(\Phi)$ will be denoted by 
$\Th_0$, and we will refer to it as the canonical
Liouville 1--form on $ST^*\M$.

   From \eq\meight\ it follows that (see \cite{\Hector})
$$
\Th_0 = \sum_{i=1}^m (p^i +\pi p^i)\, dq^i + 
\sum_{\a=1}^n (\eta^\a +\pi \eta^\a) \,d\th^\a.
\eqlabel\mnine
$$

  On the other hand, if $\Psi\:T^*\M\to ST^*\M$ is the
canonical closed imbedding of $T^*\M$ which locally is
given by
$$
\matrix{
q^i &\mapsto & q^i  \qquad\qquad  
&\th^\a &\mapsto & \th^\a,   \cr      
p^i &\mapsto &p^i  \qquad\qquad  
&\eta^\a &\mapsto  &\eta^\a,   \cr
\pi\eta^\a &\mapsto &0  \qquad\qquad 
&\pi p^i &\mapsto &0,   \cr}
$$
then the restriction of $\Th_0$ to $T^*\M$, that
will also be denoted by $\Th_0$, is locally 
given by
$$
\Th_0 = \sum_{i=1}^m p^i\, dq^i + 
\sum_{\a=1}^n \eta^\a \,d\th^\a,
$$
where, to be precise $dq^i$ and $d\th^\a$ stand for
$d\bigl(\phi^*(q^i)\bigr)$ and 
$d\bigl(\phi^*(\th^\a)\bigr)$, respectively.

  The canonical Liouville 1--form was defined in 
\cite{\Sanchez}\ in a different way, which is equivalent
to ours: 
 
\proclaim Teorema 4.1.
The canonical Liouville 1--form $\Th_0$ is the
only $\Pi$--semibasic 1--form on $ST^*\M$
that satisfy
$$
\Si_\om^*(\Th_0)= \om  \qquad\forall\om\in\O^1(\A),
$$
where $\Si_\om$ is the section of the cotangent
superbundle corresponding to $\om$.

\Proof:
It is enough to work on a local chart of $\M$. Thus,
if $\om=\sum_{i=1}^m w^i \,dq^i + 
\sum_{\a=1}^n \om^\a \,d\th^\a$ on an open subset $\U$
of $M$, we have
$$
\leqalignno{  
\Si_\om^*(\Th_0)
&= \si_\om^\sharp(\widehat\Th_0)  \cr
&= \si_\om^\sharp\Bigl(
\sum_{i=1}^m \si_\om^*(p^i +\pi p^i)\, d\hat q^i + 
\sum_{\a=1}^n\si_\om^*(\eta^\a+\pi\eta^\a)\,d\hat\th^\a 
\Bigr)  \cr
&= \si_\om^\sharp\Bigl(
\sum_{i=1}^m w^i \,d\hat q^i + 
\sum_{\a=1}^n \om^\a \,d\hat\th^\a \Bigr)  \cr
&= \sum_{i=1}^m w^i \,d\bigl(\si_\om^*(q^i)\bigr) + 
\sum_{\a=1}^n \om^\a \,d\bigl(\si_\om^*(\th^\a) \bigr)
\cr
&= \sum_{i=1}^m w^i \,dq^i + 
\sum_{\a=1}^n \om^\a \,d\th^\a = \om.  \cr}
$$

  On the other hand, a general graded 1--form $\Th$
on $ST^*\M$ is written locally as
$$
\Th = \sum_{i=1}^m A^i \,dq^i +\sum_{i=1}^m B^i\,dp^i
\sum_{\a=1}^n C^\a\,d\pi\eta^\a 
+\sum_{\a=1}^n D^\a\,d\th^\a
+\sum_{\a=1}^n E^\a\,d\eta^\a
+\sum_{i=1}^m F^i\,d\pi p^i,
$$
but, if it is $\Pi$--semibasic then
$B^i$, $C^\a$, $E^\a$ and $F^i$ vanish, and the
previous argument fix the other two supercoordinates,
and the uniqueness follows.   \qed

\declare Remark 4.1.
 Although $\Th_0$ is formally equal to the canonical 
1--form of the cotangent bundle in non--graded geometry, 
it turns out that the graded 2--form $-d\Th_0$ is 
degenerate; nevertheless, if one restricts $\Th_0$ to the 
cotangent supermanifold $T^*\M$, then
$-d\Th_0$ is a non--degenerate graded 2--form that 
will be called the canonical graded 2--form and 
will be denoted by $\O_0$.  We
refer to~\cite{\Sanchez} for details.

\beginsection 5. The super Legendre transformation.

\subsection 5.1. The vertical superendomorphism.

   As in the non--graded case,  in order to define
intrinsically the vertical superendomorphism, we
need to define vertical lifts.  We shall accomplish
this generalizing the ideas of the non--graded
case (see, for instance \cite{\Yano}).

   Let $\U$ be an open subset of $M$ such that
$\bigl(\U,\AU\bigr)$ is isomorphic to a superdomain.
We associate to each superfunction $f\in\AU$ the
superfunction $f^V\in T\AU$ defined by
$$
f^V := \sum^m_{i=1} \pd{F}{q^i}v^i 
+ \sum^n_{\a=1} \pd{F}{\th^\a}\z^\a, 
$$
where $F:= \tau^*(f)\in T\AU$.  It turns out that,
any supervector field $Y$ on $T\M$ is determined
by its action on the superfunctions $f^V$: 

\proclaim Lema 5.1.
If $Y\in\X(T\A)$ satisfy
$$
Y(f^V) = 0  \qquad \forall f\in\AU,
$$
then $Y\equiv 0$ on $\tau^{-1}(\U)$.

\Proof:
If the local expression for $Y$ is
$$
Y = \sum_{k=1}^m A^k\partial_{q^k} +
\sum_{k=1}^m B^k\partial_{v^k} 
+ \sum_{\ga=1}^n C^\ga\partial_{\th^\ga}
+ \sum_{\ga=1}^n D^\ga\partial_{\z^\ga},
$$
then
$$
\leqalignno{
0 = Y(f^V) &= 
\sum_{k,i}^m A^k\sopd{F}{q^k}{q^i}v^i
+\sum_{k,\a}^m A^k
\sopd{F}{q^k}{\th^\a}\z^\a
+\sum_{k=1}^m B^k\pd{F}{v^k}
& (\ecnum) \cr
&\quad +\sum_{\ga,i}^m C^\ga
\sopd{F}{\th^\ga}{q^i}q^i
+ \sum_{\ga,\a}^m C^\ga
\sopd{F}{\th^\ga}{\th^\a}\z^\a
+ \sum_{\ga=1}^n D^\ga\pd{F}{\z^\ga}.  \cr}
\eqlabel\mten
$$
Plugging $f=q^j$ in \eq\mten, one 
gets $B^j=0$; similarly, if 
$f=\th^\b$ it follows that $D^\b=0$.

On the other hand, taking $f=q^lq^j$
in \eq\mten, one gets
$$
A^l v^j + A^j v^l=0;
$$
in particular, if $l=j$, then 
$A^j v^j=0$, and therefore $A^j=0$.
Similarly, using $f=q^j\th^\b$
one gets $C^\b=0$, and the lemma is proved. \qed

\declare Definition 5.1.
If $X$ is a supervector field on $\M$, its 
vertical lift is the supervector field  
$X^V$ on $T\M$ defined by
$$
X^V(f^V) = \tau^*\bigr(X(f)\bigl) 
\qquad \forall f\in\A.
$$

Similarly, if $X$ is a supervector field along $\Tau$,
then we define its vertical lift by the relations
$$
X^V(f^V) = X(f) 
\qquad \forall f\in\AU.
\eqlabel\mtwelve
$$

  In local supercoordinates, if 
$X = \sum_{i=1}^m X^i\partial_{\hat q^i} 
+ \sum_{\a=1}^n \chi^\a\partial_{\hat\th^\a}$, then
$$
X^V = \sum_{i=1}^m X^i\partial_{v^i} 
+ \sum_{\a=1}^n \chi^\a\partial_{\z^\a}.
$$

  The situation is slightly different in the
tangent superbundle $ST\M$.  The natural thing to
do is to replace  $f^V$ by the superfunctions
$$
f^V := 
\sum^m_{i=1} \pd{F}{q^i}\, (v^i + \pi v^i) 
+ \sum^n_{\a=1} \pd{F}{\th^\a}\,(\z^\a + \pi\z^\a).
\eqlabel\meleven 
$$
Eventhough a general supervector
field is not determined by its action on these
superfunctions (for instance,  $Y(f^V) =0$ por all 
$f\in\AU$ if $Y=\partial_{v^i}-\partial_{\pi v^i})$, 
one can check, using the same argument as before, 
that homegenous supervector fields are determined by 
its action on superfunctions of the form \eq\meleven.

  Thus we define the vertical lift of an homogeneous
supervector field $X\in\X(\A)$ as the supervector field 
$X^V\in\X(ST\A)$ that satisfies
$$
X^V(f^V) = \tau^*_0\bigr(X(f)\bigl), 
\qquad \forall f\in\A.
$$
Moreover, if $X=X_0+X_1$, then we define
$X^V:=X^V_0+X^V_1$.

  Similarly, if $X$ an homogeneous supervector field
along the canonical projection of $ST\M$ onto $\M$,
its vertical lift is also defined by the equation 
\eq\mtwelve, where now $\Tau$ denotes the projection
of $ST\M$, and, of course, in the general case by
$X^V:=X^V_0+X^V_1$. 

 We are now in a position to define, in an intrinsic
way, the two objects that encodes all the geometric
information of the tangent superbundle.

\declare Definition 5.2.
The vertical superendomorphism is
the graded tensor field of type $(1,1)$ 
$S \: \X(ST\A) \to \X(ST\A)$  defined by
$$
S(Y) := T\tau(Y)^V.
\eqlabel\mthirdt
$$

  The morphism of
$T\A$--modulos $S \: \X(T\M) \to \X(T\M)$, defined also
by  \mthirdt, except for, now $\Tau$ denotes the
restriction to $T\M$ will also be called vertical 
superendomorphism.

  On the other hand, if 
$$
\leqalignno{
Y &= \sum_{i=1}^m Y^i\partial_{q^i}
+ \sum_{i=1}^m{\cal Y}^i\partial_{v^i}
+ \sum_{\a=1}^n \tilde\Xi^\a\partial_{\pi\z^\a} 
& (\ecnum)  \cr
&\quad
+ \sum_{\a=1}^n \Upsilon^\a\partial_{\th^\a} 
+ \sum_{\a=1}^n \Xi^\a\partial_{\z^\a}  
+ \sum_{i=1}^m\tilde{\cal Y}^i\partial_{\pi v^i},
\cr}
$$ 
then, using the change rule~\cite{\Leites},
$$
SY = \sum_{i=1}^m Y^i\partial_{v^i} 
+ \sum_{\a=1}^n \Upsilon^\a\partial_{\z^\a} 
+ \sum_{i=1}^m Y^i\partial_{\pi v^i} 
+ \sum_{\a=1}^n \Upsilon^\a\partial_{\pi\th^\a}.
\eqlabel\mfourt
$$
In particular, it is clear that 
$$
\Im S=\ker S
=\set{Y\: Y\,
\hbox{is vertical with respect to}\,\Tau_0},
$$
and that the matrix of $S$, in terms of the 
supercoordinates we have been using, is
$$
S = \pmatrix{
0 & 0 & 0 & 0 & 0 & 0  \cr
I & 0 & 0 & 0 & 0 & 0  \cr
0 & 0 & 0 & I & 0 & 0  \cr
0 & 0 & 0 & 0 & 0 & 0  \cr
0 & 0 & 0 & I & 0 & 0  \cr
I & 0 & 0 & 0 & 0 & 0  \cr},
$$
while the corresponding matrix for the vertical
superendomorphism of $T\M$ would be
$$
S = \pmatrix{
0 & 0  & 0 & 0   \cr
I & 0  & 0 & 0   \cr
0 & 0  & 0 & 0   \cr
0 & 0  & I & 0   \cr}.
$$

\declare Definition 5.3.
The Liouville supervector field $\De$ is the vertical
lift of the total time derivative. 
In other words, $\De$ is the supervector field on
$\X(ST\A)$ (or $\X(T\A)$) defined by
$$
\De= \T^V.
$$

\subsection 5.2. Graded Cartan forms.

  In analogy with ordinary Lagrangian mechanics, the Cartan 
graded 1--form associated to a given Lagrangian superfunction
$L$ in $ST\A$ is defined by 
$$
\Th_L := dL \circ S.
$$
   Using \eq\mfourt\ it is easy to check that in local 
supercoordinates
$$
\Th_L = 
\Bigr(\pd{L}{v^i} -(-1)^{\abs{L}} 
\pd{L}{\pi v^i}\Bigl) dq^i +
\Bigr(\pd{L}{\pi\z^\a} -(-1)^{\abs{L}} 
\pd{L}{\z^\a} \Bigl) d\th^\a.
$$
  The cartan graded 2--form is defined as 
the exact graded 2--form
$$
\O_L = -d\Th_L,
$$
hence in local supercoordinates is written as
$$
\leqalignno{
-\O_L = 
& \Bigl(\sopd{L}{q^i}{v^j}
- (-1)^{\abs{L}}\sopd{L}{q^i}{\pi v^j}\Bigr)
dq^i \wedge dq^j \cr
&+ \Bigl(\sopd{L}{v^i}{v^j} 
- (-1)^{\abs{L}}\sopd{L}{v^i}{\pi v^j}\Bigr)
dv^i \wedge dq^j \cr  
&+ \Bigl(\sopd{L}{\pi\z^\a}{v^j} 
- (-1)^{\abs{L}}\sopd{L}{\pi\z^\a}{\pi v^j}\Bigr)
d\pi\z^\a \wedge dq^j \cr
&- \Bigl((-1)^{\abs{L}}\sopd{L}{\th^\a}{v^j} 
+ \sopd{L}{\th^\a}{\pi v^j}\Bigr)
d\th^\a \wedge dq^j  \cr
&- \Bigl((-1)^{\abs{L}}\sopd{L}{\z^\a}{v^j} 
+ \sopd{L}{\z^\a}{\pi v^j}\Bigr)
d\z^\a \wedge dq^j  \cr
&- \Bigl((-1)^{\abs{L}}\sopd{L}{\pi v^i}{v^j} 
+ \sopd{L}{\pi v^i}{\pi v^j}\Bigr)
d\pi v^i \wedge dq^j  \cr
&+ \Bigl(\sopd{L}{q^i}{\pi\z^\b} 
- (-1)^{\abs{L}}\sopd{L}{q^i}{\z^\b}\Bigr)
dq^i \wedge d\th^\b
&  (\ecnum)\cr
&+ \Bigl(\sopd{L}{v^i}{\pi\z^\b} 
- (-1)^{\abs{L}}\sopd{L}{v^i}{\z^\b}\Bigr)
dv^i \wedge d\th^\b  \cr
&+ \Bigl(\sopd{L}{\pi\z^\a}{\pi\z^\b} 
- (-1)^{\abs{L}}\sopd{L}{\pi\z^\a}{\z^\b}\Bigr)
d\pi\z^\a \wedge d\th^\b   \cr
&- \Bigl((-1)^{\abs{L}}\sopd{L}{\th^\a}{\pi\z^\b} 
+ \sopd{L}{\th^\a}{\z^\b}\Bigr)
d\th^\a \wedge d\th^\b  \cr
&- \Bigl((-1)^{\abs{L}}\sopd{L}{\z^\a}{\pi\z^\b} 
+ \sopd{L}{\z^\a}{\z^\b}\Bigr)
d\z^\a \wedge d\th^\b  \cr
&- \Bigl((-1)^{\abs{L}}\sopd{L}{\pi v^i}{\pi\z^\b} 
+ \sopd{L}{\pi v^i}{\z^\b}\Bigr)
d\pi v^i \wedge d\th^\b .  \cr}
$$
Therefore, the matrix associated to $\O_L$
is of the form
$$
\O_L = \pmatrix{
A_1   & A_2 & A_3 & B_1   & B_2 & B_3  \cr
-A^t_2& 0   & 0   & B_4   & 0   & 0    \cr
-A^t_3& 0   & 0   & B_ 5  & 0   & 0    \cr
C_1   & C_4 & C_5 & D_1   & D_2 & D_3  \cr
C_2   & 0   & 0   & D^t_2 & 0   & 0    \cr
C_3   & 0   & 0   & D^t_3 & 0   & 0    \cr},
$$
where $C_i = -(-1)^{\abs{L}}B^t_i$; in particular, 
$\O_L$ will be degenerate for every 
superfunction $L\in ST\M$.

\subsection 5.3. The super--Legendre transformation.

   If $Y$ is a vertical supervector field with respect to 
$\Tau$ (i.e. $Y\circ \tau^* = 0$) then $\Th_L(Y) =0$,  
and therefore $\Th_L$ is a
$\Tau$--semibasic graded 1--form, and since $\Tau$ is a
submersion, it has associated a unique graded 1--form 
$\widehat{\Th}_L$ along $\Tau$~\cite{\Hector}.  In terms 
of the basis $\set{d\hat q^i, d\hat \th^\a}$, 
$\widehat{\Th}_L$ has the same coordinates as $\Th_L$
corresponding to the elements $\set{dq^i, d\th^\a}$
(which is not a full basis of $\O^1\AU$), hence
$$
\widehat{\Th}_L = 
\Bigr(\pd{L}{v^i} -(-1)^{\abs{L}} \pd{L}{\pi v^i}\Bigl) 
d\hat q^i +
\Bigr(\pd{L}{\pi\z^\a} -(-1)^{\abs{L}} \pd{L}{\z^\a} \Bigl) 
d\hat\th^\a.
$$

   In analogy with non--graded geometry, see
\cite{\Pepinii}, the section $\FL \: ST\M \to ST^*\M$ 
along $\Tau$ that corresponds to the graded 1--form 
$\widehat{\Th}_L$ could be considered as the Legendre 
transformation, but in view of the degeneracy of
$\O_L$ for every $L\in ST\M$, we shall restrict our 
attention to the case when the super--Lagrangian 
$L\in T\M\subset ST\M$,
(i.e. when $L$ does not depend on the variables $\pi v^i$
or $\pi\th^\a$) and consider the restriction of $\FL$ to 
$T\M$.

\declare Definition 5.4.
If $L$ is a super--Lagrangian in $T\M$, the super--Legendre 
transformation associated to $L$ is the restriction of the
map $\FL$ to $T\M$.  We shall denote the super Legendre 
transformation by $FL$. Hence
$$
FL \: T\M \longrightarrow ST^*\M.
$$

   When $L\in T\M$ the matrix of $\O_L$ reduces to
$$
\O_L = \pmatrix{
A_1   & A_2 & B_1   & B_2  \cr
-A^t_2& 0   & B_4   & 0     \cr
C_1   & C_4 & D_1   & D_2   \cr
C_2   & 0   & D^t_2 & 0     \cr},
$$
and to analyze its degeneracy it is necessary to consider 
the parity of $L$. If $L$ is even then $\O_L$ is 
non--degenerate if, and only if, the matrices $A_2$ and 
$D_2$ are invertible; in other words, exactly when
$$
\sopd{L}{v^i}{v^j} \sepword{and}\quad 
\sopd{L}{\z^\a}{\z^\b} 
\sepword{are invertible.} 
\eqlabel\THREE
$$

   We also notice that if $\abs{L} = 0$ then $FL$ takes
values in $T^*\M$.  In fact, locally $FL = (fl,fl^*)$ 
is determined by the morphism of superalgebras 
$fl^*\: T^*\AU \to T\AU$ described by the relations:
$$
\matrix{
q^i &\mapsto &q^i  \qquad\qquad  &\th^\a &\mapsto &\th^\a,  
 \cr      
p^i &\mapsto &\pd{L}{v^i}  \qquad\qquad  
&\eta^\a &\mapsto &-\pd{L}{\z^\a},    \cr}
$$
which, by the inverse function theorem~\cite{\Leites},
will be a local diffeomorphism when the Jacobian is
invertible, and this happens exactly when \eq\THREE\
holds.

   On the other hand, if $L$ is odd, $\O_L$ is
non--degenerate if, and only if, the off diagonal terms
are non--degenerate. This implies that $m=n$ and that
$B_2$ is invertible.  In other words, that
$$
\sopd{L}{\z^\a}{v^j}\sepword{is invertible.} 
\eqlabel\FOUR
$$

   Unlike the even case, the super--Legendre 
transformation does not take values in $T^*\M$, but on
the subsupermanifold of $ST^*\M$ of dimension $(m+n,n+m)$
obtained by imposing the conditions
$$
p^i = 0 \quad 1\leq i\leq m \sepword{and}\quad
\eta^\a = 0 \quad 1\leq \a\leq n.
$$
Moreover, locally $FL$ is given by the assignments
$$
\matrix{
q^i &\mapsto &q^i  \qquad\qquad  &\th^\a &\mapsto &\th^\a,  
 \cr      
\pi\eta^\a &\mapsto &\pd{L}{\z^\a}  \qquad\qquad  
&\pi p^i &\mapsto &-\pd{L}{v^i};    \cr}
$$
nevertheless, when $m=n$, again by the inverse function
theorem, $FL$ is a local diffeomorphism exactly when
\eq\FOUR\ holds.  We have, therefore, proved the following

\proclaim Proposition 5.1.
The super--Legendre transformation $FL$ is a local 
diffeomorphism if, and only if, the graded form $\O_L$
is non--degenerate. In either case, we say that the super
Lagrangian $L$ is regular.

  The super--Legendre transformation has the same properties 
as the usual Legendre transformation~\cite{\Abraham}.

\proclaim Proposition 5.2.
Let $L$ be a super--Lagrangian in $ST\A$, then 
$\FL^*(\Th_0) = \Th_L$.  Moreover, when $L\in T\A$ and
one restricts $\Th_L$ and $\Th_0$ to the appropriate
subsupermanifolds (for instance to $T\M$ and $T^*\M$
respectively, when $\abs{L} =0$) then also
$FL^*(\Th_0) = \Th_L$.

\Proof:
  This is immediate from the local coordinate expressions.
Let us simply remind the reader that, for instance, when
$\abs{L} =0$ then $\FL = (fl,fl^*)$ is the morphism of 
supermanifolds associated to the morphism of superalgebras 
$fl^*\:ST^*\AU \to  ST\AU$ given by
$$
\matrix{
q^i &\mapsto &q^i  \qquad\qquad  &\th^\a &\mapsto &\th^\a,  
 \cr      
p^i &\mapsto &\pd{L}{v^i}  \qquad\qquad  
&\eta^\a &\mapsto &-\pd{L}{\z^\a},    \cr
\pi\eta^\a &\mapsto &\pd{L}{\pi\z^\a}   \qquad\qquad
&\pi p^i &\mapsto &-\pd{L}{\pi v^i}.    \cr}
$$  \qed

 When $L$ is a regular super--Lagrangian there exists a unique
supervector field $\Ga_L$ in $\X(\M)$ such that 
$$
i_{\Ga_L}\O_L = d\,E_L,
$$
where the superenergy is defined by $E_L := \De L - L$ and
$\De$ is the Liouville supervector field.  Moreover, $\Ga_L$
is a super Second Order Differential Equation,
see~\cite{\Ibort} for details.

\proclaim Proposition 5.3.
Let $L$ be a super--Lagrangian in $T\A$ such that $FL$ is
a diffeomorphism (in such case we say $L$ is hyperregular).
Then $V= (FL^{-1})^* \circ \Ga_L \circ FL^*$ is a Hamiltonian
supervector field with Hamiltonian $H:= (FL^{-1})^* E_L$.
Reciprocally, if $H$ is the superfunction $H:= (FL^{-1})^* E_L$,
then the Hamiltonian supervector field $V$ associated to $L$ is
$FL$--related to $\Ga_L$.

\Proof:
Let $X$ be a supervector field on $T^*\M$. Since $FL$ is
a diffeomorphism there exists $Y\in \X(T\A)$ such that 
$X= (FL^{-1})^* \circ Y \circ FL^*$. Using Lemma 3.1 twice
we have
$$
\leqalignno{
i_V\O_0(X) = \O_0(V,X) 
&= (FL^{-1})^*\bigr[ FL^*(\O_0)(\Ga_L,Y)\bigr] 
& (\ecnum) \cr      
&= (FL^{-1})^*\bigr[ i_{\Ga_L}\O_L(Y)\bigr]  \cr
&= \bigr[d (FL^{-1})^* (E_L)\bigr] (X) = d\,H(X),  \cr}
$$
and the first assertion follows.

  As for the second statement, we consider
$Z= (FL^{-1})^* \circ \Ga_L \circ FL^*$; then the previous
argument gives that $i_Z\O_0 =d\,H$, and since $\O_0$ is
non--degenerate, $Z= V$, and the proposition is proved.  \qed

  Moreover, since $\Ga_L$ is a super SODE then
$i_{\Ga_L}\Th_L = \De(L) =: A$, and the same argument gives
us $\Th_0(V) = (FL^{-1})^*A$ when $L$ is hyperregular.  The
correspondence between the Lagrangian and Hamiltonian
formulations in supermechanics is clear.

\bigskip 

\noindent {\bf Acknowledgements} 
\smallskip
 
We thank J. Monterde and O.A. S\'anchez--Valenzuela 
for useful comments related to the tangent superbundle.
Partial finantial support from DGICYT under projects  
PS--90.0118 and PB--93.0582 is acknowledged. 
HF thanks the Vicerrector\1a de 
Investigaci\'on de la Universidad de Costa Rica and 
the Agencia Espa\~nola de Cooperaci\'on Internacional.
Finally we tank the first referee for the suggestion of
references 10 and 11.
  
\bigskip  
\noindent {\bf References} 
\frenchspacing
\bigskip

\refno\Abraham.
R. Abraham and J.E. Marsden, {\it Foundations of Mechanics}, 
Second Edition. The Benjamin Cumming Publishing Company, Reading,
Massachusetts (1978).

\refno\Batchelor.
M. Batchelor, {\it The structure of supermanifolds}, Trans. 
Am. Math. Soc.\ {\bf253} (1979), 329--338.

\refno\Sancheziii.
C.P. Boyer and O.A. S\'anchez--Valenzuela, {\it Lie supergroup actions on 
 supermanifolds}, Trans. Amer. Math. Soc. {\bf 323} (1991) 151--175.
 
 \refno\Hector.
J.F. Cari\~nena and H. Figueroa, {\it A geometrical version 
 of Noether's theorem in supermechanics}. Reports on 
 Mathematical Physics  {\bf 34} (1994) 277--303.

 \refno\Hectoriii.
J.F. Cari\~nena and H. Figueroa, {\it Higher order
Lagrangian mechanics in the graded context}, preprint DFTUZ
(1996).

\refno\Pepino.
J.F. Cari\~nena, C. L\'opez and E. Mart\1nez, 
{\it A new approach to the converse of Noe\-ther's theorem},
 J. Phys. A: Math. Gen. {\bf 22} (1989) 77--86.
 
 \refno\Pepinii.
J.F. Cari\~nena, C. L\'opez and E. Mart\1nez, 
{\it Sections along a map applied to Higher Order Lagrangian Mechanics. 
Noether's theorem}, Acta Applicandae Mathematicae\ {\bf25} (1991) 127--151.

\refno\Pepiniv.
J.F. Cari\~nena, J. Fern\'andez--N\'u\~nez,  
{\it Geometric theory of time--dependent singular Lagrangians}, 
Fortschr. Phys. {\bf41:6} (1993) 517--552.

\refno\Gawedzki.
K. Gawedzki, {\it Supersymmetries--mathematics of supergeometry},
Ann. Inst. Henri Poin\-ca\-r\'e, vol XXVII, {\bf4} (1977) 335--366.

\refno\Hilli.
D. Hill and S.R. Simanca, {\it Newlander--Niremberg theorem on supermanifolds
with boundary\/}, Rivista di Matematica, Universit\`a di Parma, (5) {\bf 2}
(1993) 213--228

\refno\Hillii.
D. Hill and S.R. Simanca, {\it The supercomplex Frobenius theorem},
 Annales Polinici Mathematici,  {\bf 55}
(1991) 139--155

\refno\Ibort.
L.A. Ibort and J. Mar\1n--Solano, {\it Geometrical foundations of 
Lagrangian supermechanics and supersymmetry}, Reports on Mathematical 
Physics {\bf 32:3} (1993) 385--409.

\refno\Kostant.
B. Kostant, {\it Graded manifolds, graded Lie theory and
prequantization}, in  Differential Geometrical Methods in 
Mathematical Physics, Lecture Notes in Mathematics {\bf570},
Springer, Berlin, 1977.

\refno\Leites.
D.A. Le\u{\i}tes, {\it Introduction to the theory of 
supermanifolds}, Russian Math. Surveys\ {\bf35:1} (1980) 1--64.

\refno\Manin.
Yu.I. Manin, {\it Gauge Field theory and Complex Geometry},
Nauka. Moscow, (1984). English transl., Springer--Verlag, New York, (1988).

\refno\MonterdeSanchez.
J. Monterde y O. A. S\'anchez--Valenzuela, {\it Existence
and uniqueness of solutions to superdifferential equations},
J. Geom. and Phys. {\bf 10} (1993) 315--343.

\refno\Monterde.
J. Monterde y O. A. S\'anchez--Valenzuela, {\it On 
the Batchelor trivialization of the tangent 
supermanifold} en Proccedings of the fall Workshop on
Differential Geometry and its Applications, Barcelona
(1993) 9--14.

\refno\Sanchez.
O.A. S\'anchez--Valenzuela, {\it Differential--Geometric approach
to  supervector bundles}, Comunicaciones T\'ecnicas IIMASS--UNAM, 
(Serie Naranja) {\bf457}, M\'exico, (1986).

\refno\Sanchezi.
O.A. S\'anchez--Valenzuela, 
{\it On Grassmanian supermanifolds}, Trans. Amer. Math. Soc. {\bf307} (1988)
597--614.

\refno\Sanchezii.
O.A. S\'anchez--Valenzuela, {\it Linear supergroup actions I;
on the defining properties}, Trans. Amer. Math. Soc. {\bf307} (1988)
569--595.

\refno\Tennison.
B.R. Tennison, {\it Sheaf
 theory}, London Math. Soc. Lecture Note Series {\bf20}, Cambridge Univ. 
 Press, Cambridge, (1975).
 
 \refno\Yano.
K. Yano and S. Ishihara, {\it Tangent and Cotangent bundles}, 
 Marcel Decker, Inc. New York, (1973).

\bye